\documentclass[11pt,a4paper]{article}
\usepackage[T1]{fontenc}
\usepackage[utf8]{inputenc}
\usepackage[english]{babel}
\usepackage{lmodern}
\usepackage[margin=28mm]{geometry}
\usepackage{mathtools,amssymb,amsthm}
\usepackage{enumitem}
\usepackage{microtype}
\usepackage[hidelinks]{hyperref}

\newtheorem{theorem}{Theorem}
\newtheorem{lemma}[theorem]{Lemma}

\theoremstyle{remark}
\newtheorem{remark}[theorem]{Remark}

\DeclareMathOperator{\dom}{dom}
\DeclareMathOperator{\Graph}{Graph}
\DeclareMathOperator{\Epi}{Epi}
\DeclareMathOperator{\cl}{cl}
\DeclareMathOperator{\conv}{conv}
\DeclareMathOperator{\Supp}{Supp}
\DeclareMathOperator{\spn}{span}

\newcommand{\R}{\mathbb{R}}
\newcommand{\N}{\mathbb{N}}
\newcommand{\pair}[2]{\langle #1,#2\rangle}
\setlist[enumerate]{leftmargin=*,itemsep=0.4em,topsep=0.5em}
\hypersetup{
  pdftitle={Maximal convex subdifferential spaces},
  pdfauthor={M. D. Voisei},
  pdfsubject={Subdifferential maximality in locally convex spaces},
  pdflang={en-US}
}

\begin{document}
\title{Spaces with maximal monotone convex
	subdifferentials}
\author{M.\,D. Voisei}
\date{}
\maketitle

\begin{abstract}
We characterize the real Hausdorff locally convex spaces on which every
proper lower semicontinuous convex function has a maximal monotone convex 
subdifferential. Equivalent conditions are given in terms of the portability
of closed convex sets and identities involving upper envelopes. Among normed
spaces, this property is equivalent to completeness; among Fr\'echet spaces,
it is equivalent to normability. We also prove that the property is preserved
under products with finite-dimensional spaces.
\end{abstract}



\section{Introduction and preliminaries}

The aim of this paper is to characterize the Hausdorff locally convex spaces
on which every proper lower semicontinuous convex function has a maximal
monotone convex subdifferential. A real Hausdorff locally convex space $(X,\tau)$,
with continuous dual $X^*$, has the \emph{subdifferential maximality property}
(SM) if $\partial f:X\rightrightarrows X^*$ is maximal monotone for every
$f\in\Gamma_\tau(X)$. Here $\Gamma_\tau(X)$ denotes the class of proper
$\tau$-lower semicontinuous convex functions $f:X\to\overline{\R}$,
and monotonicity is taken with respect to the canonical pairing
$\pair{x}{x^*}=x^*(x)$.

The problem of characterizing spaces with (SM) was posed as Problem (P1)
in \cite{mmsd}. We extend the characterizations in \cite[Theorem~17]{mmsd}
by proving that (SM) is equivalent to the portability of every nonempty
closed convex subset of $X$. We also obtain characterizations through
upper-envelope identities, both on $X$ and on the effective domains of the
functions; see Theorems~\ref{ch} and~\ref{ch-dom}. We then characterize
(SM) among normed and Fr\'echet spaces and establish its invariance under
products with finite-dimensional spaces.

We use $\overline{\R}=\R\cup\{-\infty,+\infty\}$ and the conventions
$\sup\varnothing=-\infty$ and $\inf\varnothing=+\infty$.
For a proper function $f$, set
\[
 \dom f=\{x\in X\mid f(x)<+\infty\},\qquad
 \Epi f=\{(x,t)\in X\times\R\mid f(x)\le t\}.
\]
For an operator $T:X\rightrightarrows X^*$, $\Graph T$ denotes its graph,
$D(T)=\{x\mid T(x)\ne\varnothing\}$ its domain, and
$R(T)=\bigcup_{x\in X}T(x)$ its range. The indicator $\iota_C$ equals
$0$ on $C$ and $+\infty$ outside $C$.

\subsection{Compatible topologies and duality}

Property (SM) depends only on the duality $(X,X^*)$. Indeed, if $\mu$ is
a locally convex topology compatible with this duality, then
$\Gamma_\mu(X)=\Gamma_\tau(X)$: closed convex sets, and hence closed
convex epigraphs, are the same for compatible locally convex topologies.
The corresponding subdifferentials also coincide. In particular, one may
replace $\tau$ by the weak topology $\sigma(X,X^*)$.

More generally, let $(X,Y,c)$ be a real dual system, where
$c:X\times Y\to\R$ is a bilinear pairing. Write $\Gamma(X,Y,c)$ for the
proper convex $\sigma(X,Y)$-lower semicontinuous functions. For such a
function, define
\begin{equation}\label{eq:subdifferential}
 \partial_c f(x)=
 \{y\in Y\mid f(z)\ge f(x)+c(z-x,y)\text{ for all }z\in X\},
 \qquad x\in\dom f,
\end{equation}
and $\partial_c f(x)=\varnothing$ otherwise. The dual system has (SM) if
every such subdifferential is maximal monotone with respect to $c$.

This property is symmetric: $(X,Y,c)$ has (SM) if and only if $(Y,X,c')$
does, where $c'(y,x)=c(x,y)$. To see this, define
\[
 f^*(y)=\sup_{x\in\dom f}\{c(x,y)-f(x)\}.
\]
Fenchel--Moreau duality makes conjugation a bijection between
$\Gamma(X,Y,c)$ and $\Gamma(Y,X,c')$, with inverse conjugation, and
Fenchel--Young equality gives
\begin{equation}\label{eq:inverse-subdifferential}
 \partial_{c'}f^*=(\partial_c f)^{-1}.
\end{equation}
Maximal monotonicity is preserved under inversion. These assertions also
hold for nonseparated dual systems: pass to the separated quotients by
$X_0=\{x\mid c(x,Y)=0\}$ and $Y_0=\{y\mid c(X,y)=0\}$, observing that
weakly lower semicontinuous functions are constant on the corresponding
kernel cosets.

From now on, unless explicitly stated otherwise, $(X,\tau)$ is a real
Hausdorff locally convex space, and $X^*$ is equipped with
$w^*=\sigma(X^*,X)$. Conjugates and subdifferentials on $X^*$ are taken
with respect to $(X^*,X)$, not a larger dual. We omit the pairing subscript.

The Fitzpatrick function of $T:X\rightrightarrows X^*$ is
\begin{equation}\label{eq:fitzpatrick}
 \varphi_T(x,x^*)=
 \sup_{(a,a^*)\in\Graph T}
 \bigl\{\pair{x-a}{a^*}+\pair{a}{x^*}\bigr\},
 \qquad (x,x^*)\in X\times X^*;
\end{equation}
see \cite{MR1009594}.

\subsection{Portable hulls and upper envelopes}

For $C\subset X$, its normal cone is
\[
 N_C(a)=\{a^*\in X^*\mid \pair{y-a}{a^*}\le0
                         \text{ for every }y\in C\},\quad a\in C,
\]
and $N_C(a)=\varnothing$ for $a\notin C$. Thus $N_C=\partial\iota_C$.
The \emph{portable hull} of $C$ is
\begin{equation}\label{eq:portable-hull}
 C^\#=\{x\in X\mid \pair{x-a}{a^*}\le0
                     \text{ for all }(a,a^*)\in\Graph N_C\}.
\end{equation}
It is the intersection of the supporting closed half-spaces of $C$ whose
boundary meets $C$; see \cite{Voisei2019}. The set $C$ is \emph{portable}
if $C=C^\#$. With
$\Supp C=\{a\in C\mid N_C(a)\ne\{0\}\}$, one has
$\Supp\varnothing=\Supp X=\varnothing$,
$\varnothing^\#=X$, and
\begin{equation}\label{eq:support-points}
 x\in C^\#\quad\Longleftrightarrow\quad
 \pair{x-a}{a^*}\le0
 \quad\text{for every }a\in\Supp C\text{ and }a^*\in N_C(a).
\end{equation}
Moreover,
\[
 C\subset\cl_\tau\conv C
   \subset(\cl_\tau\conv C)^\#\subset C^\#,
 \qquad C^\#\ne\varnothing.
\]
In particular, $C^\#=X$ if $\Supp C=\varnothing$.

For $f\in\Gamma_\tau(X)$, its \emph{upper envelope} is the supremum
of the affine supporting functions supplied by $\partial f$:
\begin{equation}\label{eq:upper-envelope}
 \begin{aligned}
 f^\cup(x)
   &=\sup_{(a,a^*)\in\Graph(\partial f)}
       \{\pair{x-a}{a^*}+f(a)\}\\
   &=\sup_{a^*\in R(\partial f)}
       \{\pair{x}{a^*}-f^*(a^*)\}
     =(f^*+\iota_{R(\partial f)})^*(x).
 \end{aligned}
\end{equation}
This envelope goes back to \cite[(B), p.~606]{MR0178103}.
The \emph{portable hull} of $f$ is
\begin{equation}\label{eq:function-portable-hull}
 f^\#=f^\cup+\iota_{(\dom f)^\#}.
\end{equation}
Here addition of an indicator means extension by $+\infty$ outside the
indicated set, including when $f^\cup\equiv-\infty$. In particular,
$f^\cup\le f^\#\le f$; see \cite[Theorem~5(i)]{mmsd}.

The same definitions apply on $(X^*,w^*)$. For $g\in\Gamma_{w^*}(X^*)$,
\begin{equation}\label{eq:dual-upper-envelope}
 \begin{aligned}
 g^\cup(x^*)
   &=\sup_{(a^*,a)\in\Graph(\partial g)}
       \{\pair{a}{x^*-a^*}+g(a^*)\}\\
   &=\sup_{a\in R(\partial g)}\{\pair{a}{x^*}-g^*(a)\}.
 \end{aligned}
\end{equation}
For any $K\subset X^*$,
\[
 K^\#=\{x^*\in X^*\mid \pair{a}{x^*-a^*}\le0
                      \text{ for all }(a^*,a)\in\Graph N_K\}.
\]
If $K$ is nonempty, weak-star closed, and convex, then
$N_K=(\partial\sigma_K)^{-1}$, where $\sigma_K=\iota_K^*$, so equivalently
\begin{equation}\label{eq:dual-portable-hull}
 K^\#=\{x^*\in X^*\mid \pair{a}{x^*-a^*}\le0
                    \text{ for all }(a,a^*)\in\Graph(\partial\sigma_K)\}.
\end{equation}
Also $g^\#=g^\cup+\iota_{(\dom g)^\#}$.

For $f\in\Gamma_\tau(X)$, write
\begin{equation}\label{eq:envelope-notation}
 f^\circledast=(f^*)^\cup,\qquad
 f^{\cup*}=(f^\cup)^*,\qquad
 f^\circ=(f^\circledast)^*.
\end{equation}
Using $\partial f^*=(\partial f)^{-1}$, we obtain
\begin{equation}\label{eq:dual-envelope-expansion}
 \begin{aligned}
 f^\circledast(x^*)
   &=\sup_{(a,a^*)\in\Graph(\partial f)}
       \{\pair{a}{x^*-a^*}+f^*(a^*)\}\\
   &=\sup_{a\in D(\partial f)}\{\pair{a}{x^*}-f(a)\}
     =(f+\iota_{D(\partial f)})^*(x^*),
 \end{aligned}
\end{equation}
and
\begin{equation}\label{eq:closed-envelope}
 f^\circ=(f+\iota_{D(\partial f)})^{**}
        =\cl_\tau\conv(f+\iota_{D(\partial f)}).
\end{equation}
The last expression denotes the lower semicontinuous convex hull.
Each of $f^\cup$, $f^\circledast$, $f^{\cup*}$, and $f^\circ$ is proper
if and only if $\Graph(\partial f)\ne\varnothing$. If the graph is empty,
then $f^\cup=f^\circledast\equiv-\infty$ and
$f^{\cup*}=f^\circ\equiv+\infty$. In all cases,
\begin{equation}\label{eq:envelope-order}
 f^\cup\le f^\#\le f\le f^\circ.
\end{equation}

Finally, the \emph{low-upper envelope} of $f\in\Gamma_\tau(X)$ is
\begin{equation}\label{fsp}
 f^{\smallsmile}(x)=
 \sup\{\pair{x-a}{a^*}+f(a)\mid
       (a,a^*)\in\Graph(\partial f),\ f(a)\le f(x)\},\quad x\in X;
\end{equation}
see \cite{mmsd}. The restriction $f(a)\le f(x)$ is part of this definition.

\section{Characterizations of subdifferential maximality}

We first recall the characterization on which the proofs below rely.

\begin{theorem}[{\cite[Theorem~17]{mmsd}}]\label{max-sd-sp}
Let $(X,\tau)$ be a nontrivial real Hausdorff locally convex space.
The following statements are equivalent:
\begin{enumerate}[label=\textup{(\roman*)}]
\item For every $f\in\Gamma_\tau(X)$, $\partial f$ is maximal monotone.
\item For every $f\in\Gamma_\tau(X)$, $x\in\dom f$, and $x^*\in X^*$,
      $\varphi_{\partial f}(x,x^*)\ge\pair{x}{x^*}$.
\item For every $f\in\Gamma_\tau(X)$ and $x\in\dom f$,
      $\varphi_{\partial f}(x,0)\ge0$.
\item For every $f\in\Gamma_\tau(X)$ and $x\in\dom f$,
      $f(x)=f^{\smallsmile}(x)$.
\item For every $f\in\Gamma_\tau(X)$, $f=f^{\smallsmile}$.
\item For every $f\in\Gamma_\tau(X)$ and $x\in\dom f$, there is a net
      $\{(a_i,a_i^*)\}_{i\in I}\subset\Graph(\partial f)$ such that
      $f(a_i)\le f(x)$ for all $i$, $a_i\to x$ in $(X,\tau)$, and
      \[
       f(a_i)\longrightarrow f(x),\qquad
       \pair{x-a_i}{a_i^*}+f(a_i)\longrightarrow f(x).
      \]
\item For every $f\in\Gamma_\tau(X)$ and $x\in\dom f$, there is a net
      $\{(a_i,a_i^*)\}_{i\in I}\subset\Graph(\partial f)$ such that
      $\pair{x-a_i}{a_i^*}\to0$.
\item For every $f\in\Gamma_\tau(X)$, $x\in\dom f$, and $\varepsilon>0$,
      there is $(a_\varepsilon,a_\varepsilon^*)\in\Graph(\partial f)$
      with $\pair{x-a_\varepsilon}{a_\varepsilon^*}\ge-\varepsilon$.
\end{enumerate}
In this case, for every $f\in\Gamma_\tau(X)$,
\[
 \cl_\tau D(\partial f)=\cl_\tau\dom f,\qquad
 \cl_{w^*}R(\partial f)=\cl_{w^*}\dom f^*,
\]
and both sets are convex.
\end{theorem}

The universal quantifier over $\Gamma_\tau(X)$ is essential here; these
conditions are not asserted to be equivalent for an individual function.

\begin{theorem}\label{ch}
Let $(X,\tau)$ be a nontrivial real Hausdorff locally convex space.
The following statements are equivalent:
\begin{enumerate}[label=\textup{(\alph*)}]
\item $X$ has \textup{(SM)}.
\item Every nonempty closed convex set $C\subset X$ is portable.
\item For every $f\in\Gamma_\tau(X)$, $f=f^\cup$.
\item For every $f\in\Gamma_\tau(X)$, $f=f^\#$.
\item For every $f\in\Gamma_\tau(X)$, $f=f^\circ$.
\item For every $f\in\Gamma_\tau(X)$, $\Epi f$ is portable in
      $X\times\R$ with the product topology.
\end{enumerate}
\end{theorem}

\begin{proof}
\textup{(a)} $\Rightarrow$ \textup{(c)}, \textup{(e)}.
By \cite[Proposition~12]{mmsd}, maximal monotonicity of $\partial f$
implies $f=f^\cup$ and $f^*=f^\circledast$; conjugating the latter
identity gives $f=f^\circ$. Alternatively, $\partial f$ is maximal
cyclically monotone, and the conclusion follows from
\cite[Theorem~3.13 and Remark~3.12]{gm}.

\textup{(c)} $\Rightarrow$ \textup{(d)} follows from
\eqref{eq:envelope-order}.

\textup{(d)} $\Rightarrow$ \textup{(b)}.
For nonempty closed convex $C\subset X$, direct computation gives
\[
 (\iota_C)^\cup=(\iota_C)^\#=\iota_{C^\#};
\]
see also \cite[Lemma~4]{mmsd}. Thus (d), applied to $\iota_C$, yields
$\iota_C=\iota_{C^\#}$, and hence $C=C^\#$.

\textup{(d)} $\Leftrightarrow$ \textup{(f)} follows from
\[
 (\Epi f)^\#=\Epi(f^\#);
\]
see \cite[Theorem~5(iii)]{mmsd}.

\textup{(b)} $\Rightarrow$ \textup{(a)}.
By Theorem~\ref{max-sd-sp}(iii), it suffices to prove, for every
$f\in\Gamma_\tau(X)$ and $x\in\dom f$, that
$\varphi_{\partial f}(x,0)\ge0$. In fact, we shall find
\begin{equation}\label{eq:nonnegative-pairing}
 (a,a^*)\in\Graph(\partial f)
 \quad\text{such that}\quad \pair{x-a}{a^*}\ge0.
\end{equation}
If $x$ minimizes $f$, take $(a,a^*)=(x,0)$.
Otherwise, choose $r$ with $\inf_X f<r<f(x)$, and put $C=[f\le r]$.
Then $C$ is nonempty, closed, and convex, while $x\notin C=C^\#$.
Consequently, there are $a\in C$ and $u^*\in N_C(a)$ such that
\begin{equation}\label{eq:strict-support}
 \pair{x-a}{u^*}>0,\qquad
 f(y)\le r\ \Longrightarrow\ \pair{y-a}{u^*}\le0.
\end{equation}

Define the nonempty open convex sets
\[
 U=\bigcup_{y\in\dom f}
     (f(y)-r,\infty)\times(\pair{a-y}{u^*},\infty),
 \qquad V=(-\infty,0)^2.
\]
Convexity of $U$ follows from convexity of $f$ and affinity of
$y\mapsto\pair{a-y}{u^*}$. The sets are disjoint by
\eqref{eq:strict-support}. Separating them in $\R^2$ yields
$(\lambda,\mu)\ne(0,0)$ and $\alpha\in\R$ such that
\[
 \lambda s+\mu t\ge\alpha\ge\lambda\bar s+\mu\bar t
 \quad ((s,t)\in U,\ (\bar s,\bar t)\in V).
\]
Since $U$ is closed under increasing either coordinate, $\lambda,\mu\ge0$;
taking the supremum over $V$ then gives $\alpha\ge0$.
Passing to the lower endpoints of the intervals defining $U$, we obtain
\begin{equation}\label{eq:multiplier}
 \lambda(f(y)-r)+\mu\pair{a-y}{u^*}\ge0
 \qquad (y\in\dom f).
\end{equation}
Choose $b\in\dom f$ with $f(b)<r$. If $\mu=0$, then $\lambda>0$,
and \eqref{eq:multiplier} at $b$ is impossible. Hence $\mu>0$.
If $\lambda=0$, applying \eqref{eq:multiplier} at $y=x$ gives
\[
 0\le\mu\pair{a-x}{u^*}<0,
\]
also a contradiction. Thus $\lambda>0$. At $y=a$,
\eqref{eq:multiplier} gives $f(a)\ge r$, so $f(a)=r$ because $a\in C$.
It follows that
\[
 f(y)\ge f(a)+\pair{y-a}{(\mu/\lambda)u^*}
 \qquad (y\in\dom f).
\]
The inequality is automatic outside $\dom f$. Therefore
$a^*=(\mu/\lambda)u^*\in\partial f(a)$, and
$\pair{x-a}{a^*}>0$, proving \eqref{eq:nonnegative-pairing}.
This establishes the equivalence of (a), (b), (c), (d), and (f).

\textup{(e)} $\Rightarrow$ \textup{(a)}.
Let $g\in\Gamma_{w^*}(X^*)$ and set $f=g^*\in\Gamma_\tau(X)$.
Then $f^*=g$, and (e) gives
\[
 g^*=f=f^\circ=(g^\cup)^*.
\]
The right-hand side is proper, so $g^\cup$ cannot be identically
$-\infty$. Hence $g^\cup\in\Gamma_{w^*}(X^*)$, and conjugation yields
$g=g^\cup$. Thus condition (c) holds on $(X^*,w^*)$. Applying the
already proved implication (c) $\Rightarrow$ (a) to this space shows
that $(X^*,w^*)$ has (SM). By the symmetry of (SM) under duality,
$X$ has (SM).
\end{proof}

To obtain characterizations restricted to effective domains, we first
record the following rescaling fact.

\begin{lemma}\label{lem:truncation}
Let $f\in\Gamma_\tau(X)$, let $r\in\R$ satisfy $\inf_X f<r$, and set
$g=f+\iota_{[f\le r]}$. For every $(a,a^*)\in\Graph(\partial g)$ there
is $\theta\in(0,1]$ such that $\theta a^*\in\partial f(a)$.
\end{lemma}

\begin{proof}
Fix $(a,a^*)\in\Graph(\partial g)$ and define, on $\dom f$,
$q(y)=f(y)-f(a)-\pair{y-a}{a^*}$. The subgradient inequality for $g$
gives $q(y)\ge0$ whenever $f(y)\le r$. Therefore the nonempty open
convex set
\[
 U=\bigcup_{y\in\dom f}(f(y)-r,\infty)\times(q(y),\infty)
\]
is disjoint from $(-\infty,0)^2$. As in the preceding proof,
separation gives $\lambda,\mu\ge0$, not both zero, with
\begin{equation}\label{eq:truncation-multipliers}
 \lambda(f(y)-r)+\mu\bigl(f(y)-f(a)-\pair{y-a}{a^*}\bigr)\ge0
 \quad (y\in\dom f).
\end{equation}
Since some $b$ satisfies $f(b)<r$, the case $\mu=0$ is impossible.
At $y=a$, \eqref{eq:truncation-multipliers} gives
$\lambda(f(a)-r)\ge0$. Since $f(a)\le r$, we obtain
$\lambda(f(a)-r)=0$. Consequently,
\[
 (\lambda+\mu)(f(y)-f(a))\ge\mu\pair{y-a}{a^*}
 \quad (y\in\dom f).
\]
Taking $\theta=\mu/(\lambda+\mu)\in(0,1]$ proves the claim.
\end{proof}

\begin{theorem}\label{ch-dom}
Let $(X,\tau)$ be a nontrivial real Hausdorff locally convex space.
Each of the following statements is equivalent to \textup{(SM)}:
\begin{enumerate}[label=\textup{(\alph*')},start=3]
\item For every $f\in\Gamma_\tau(X)$ and $x\in\dom f$, $f(x)=f^\cup(x)$.
\item For every $f\in\Gamma_\tau(X)$ and $x\in\dom f$, $f(x)=f^\#(x)$.
\item For every $f\in\Gamma_\tau(X)$ and $x\in\dom f$, $f(x)=f^\circ(x)$.
\end{enumerate}
\end{theorem}

\begin{proof}
Since $\dom f\subset(\dom f)^\#$, $f^\cup=f^\#$ on $\dom f$; hence
(c') and (d') coincide. Since $f\le f^\circ$, both functions equal
$+\infty$ outside $\dom f$, so (e') is equivalent to
Theorem~\ref{ch}(e). Also, Theorem~\ref{ch}(c) implies (c').

It remains to show that (c') implies (SM). Fix $f\in\Gamma_\tau(X)$,
$x\in\dom f$, and $\varepsilon>0$, and set
\[
 r=f(x)+\varepsilon/2,\qquad
 g=f+\iota_{[f\le r]}.
\]
Then $g\in\Gamma_\tau(X)$ and $g(x)=f(x)$. By (c'),
$g^\cup(x)=g(x)$, so there is $(a,a^*)\in\Graph(\partial g)$ with
\[
 f(a)+\pair{x-a}{a^*}>f(x)-\varepsilon/2.
\]
Since $f(a)\le r$, it follows that $\pair{x-a}{a^*}>-\varepsilon$.
Lemma~\ref{lem:truncation} supplies $\theta\in(0,1]$ such that
$b^*=\theta a^*\in\partial f(a)$. If $\pair{x-a}{a^*}\ge0$, then
$\pair{x-a}{b^*}\ge0$; otherwise multiplication by $\theta\le1$ increases
the pairing. In either case,
\[
 \pair{x-a}{b^*}>-\varepsilon.
\]
Theorem~\ref{max-sd-sp}(viii) now gives (SM).
\end{proof}

\section{Normed and Fr\'echet spaces}

\begin{theorem}\label{thm:normed}
A real normed space has \textup{(SM)} if and only if it is complete.
\end{theorem}

\begin{proof}
If $X$ is Banach, (SM) is Rockafellar's maximal monotonicity theorem
\cite[Theorem~A]{MR0262827}.

Conversely, suppose $X$ is incomplete. Let $\widehat X$ be its completion,
choose $z\in\widehat X\setminus X$, and take a sequence
$(x_n)_{n\ge1}\subset X$ converging to $z$ in $\widehat X$. Set
\[
 M=\cl_X\spn\{x_n\mid n\ge1\},
\]
where the closure is taken in $X$ for its norm topology.
Then $M$ is a closed separable subspace of $X$. The sequence $(x_n)$ is
Cauchy in $M$ and has no limit there, since its limit in $\widehat X$
is $z\notin X$. Thus $M$ is incomplete.

By Fonf's theorem \cite[p.~1174]{MR1216811}, there is a nonempty bounded
closed convex set $W\subset M$ with $\cl_M\spn W=M$ and no weakly
supported points. In particular, $W$ has no support points relative to
$M^*$. Indeed, a nonzero $m^*\in M^*$ cannot vanish identically on $W$,
because $\spn W$ is dense in $M$. Its restriction to $W$ is continuous,
so attainment of its supremum on $W$ would give a weakly supported point
in Fonf's sense. Therefore
\[
 N_W^M(a)=\{0\}\qquad(a\in W).
\]
Since $M$ is closed in $X$, $W$ is also closed in $X$. For $a\in W$,
\[
 x^*\in N_W^X(a)
 \quad\Longleftrightarrow\quad
 x^*|_M\in N_W^M(a)=\{0\}.
\]
Consequently,
\[
 \Graph N_W^X=W\times M^\perp,\qquad
 M^\perp=\{x^*\in X^*\mid x^*|_M=0\}.
\]
The incomplete space $M$ is infinite-dimensional, so it is unbounded.
Since $W$ is bounded, $W\subsetneq M$. Thus $N_W^X=\partial\iota_W$
has the proper monotone extension $M\times M^\perp$, the graph of
$N_M^X$. Since $\iota_W\in\Gamma(X)$, property (SM) fails.
\end{proof}

\begin{theorem}\label{thm:frechet}
Let $(X,\tau)$ be a real Fr\'echet space, that is, a complete metrizable
locally convex Hausdorff space. The following are equivalent:
\begin{enumerate}[label=\textup{(\alph*)}]
\item $(X,\tau)$ has \textup{(SM)}.
\item $\tau$ is normable.
\item Some Banach norm on $X$ induces $\tau$.
\end{enumerate}
\end{theorem}

\begin{proof}
Completeness gives (b) $\Leftrightarrow$ (c), and
Theorem~\ref{thm:normed} gives (c) $\Rightarrow$ (a).

For the remaining implication, suppose that $X$ is nonnormable, and put
$\omega=\R^\N$, $\N=\{1,2,\ldots\}$, with its product topology.
Eidelheit's quotient theorem \cite{Eidelheit1936Zur} gives a continuous
linear surjection $Q:X\to\omega$. Since both spaces are Fr\'echet,
the open mapping theorem makes $Q$ a quotient map.

Define
\begin{equation}\label{eq:omega-example}
 F(z)=\sum_{n=1}^\infty 2^{-n}e^{z_n},\qquad z=(z_n)_{n\ge1}\in\omega.
\end{equation}
As the increasing supremum of its continuous convex partial sums,
$F$ is convex and lower semicontinuous. Also $F>0$ and $F(0)=1$, so
$F\in\Gamma(\omega)$.

Every $z^*\in\omega^*$ depends on finitely many coordinates. Indeed,
$(z^*)^{-1}((-1,1))$ contains a basic neighborhood
\[
 U=\{z\in\omega\mid |z_j|<\delta_j\text{ for }j\in J\},
\]
with $J\subset\N$ finite and $\delta_j>0$. If $y_j=0$ for all $j\in J$,
then $ty\in U$ for every $t\in\R$, so $|t z^*(y)|<1$ for all $t$.
Hence $z^*(y)=0$.

Suppose $z^*\in\partial F(z)$, where $z\in\dom F$, and let $J$ be a
finite set of coordinates on which $z^*$ depends. Choose $j\notin J$.
For $t>0$, with $e_j$ the $j$th coordinate vector, put $y=z-te_j$.
Then $\pair{y-z}{z^*}=0$, whereas
\[
 F(y)=F(z)-2^{-j}e^{z_j}(1-e^{-t})<F(z),
\]
contrary to the subgradient inequality. Thus $\partial F$ has empty graph.

Now set $f=F\circ Q$. It is convex and lower semicontinuous, and
$f(0)=F(0)=1$, so $f\in\Gamma_\tau(X)$. Suppose that
$x^*\in\partial f(x)$. For $h\in\ker Q$ and $t\in\R$,
\[
 f(x)=f(x+th)\ge f(x)+t\pair{h}{x^*},
\]
so $x^*$ vanishes on $\ker Q$. Thus it factors as $x^*=z^*\circ Q$
for a linear functional $z^*$ on $\omega$. The quotient property of $Q$
implies that $z^*$ is continuous. The subgradient inequality becomes
\[
 F(Qy)\ge F(Qx)+z^*(Qy-Qx)\qquad(y\in X).
\]
Surjectivity of $Q$ gives $z^*\in\partial F(Qx)$, a contradiction.
Hence $\partial f$ has empty graph and is not maximal monotone.
Property (SM) therefore fails on $X$.
\end{proof}

\section{Products with finite-dimensional spaces}

If two Hausdorff locally convex spaces $X$ and $Y$ have (SM), it is
natural to ask whether $X\times Y$ has (SM). This holds for Banach
spaces. The converse holds in general: if $X\times Y$ has (SM), then
so do both factors. Indeed, for nonempty closed convex $C\subset X$,
\[
 N_{C\times Y}(a,b)=N_C(a)\times\{0\},\qquad
 (C\times Y)^\#=C^\#\times Y.
\]
Theorem~\ref{ch} then applies to $X$, and similarly to $Y$.
The zero space has (SM) directly, so it causes no exception here or below.

We prove that the product question has an affirmative answer whenever
one factor is finite-dimensional.

\begin{theorem}\label{XxR}
A real Hausdorff locally convex space $X$ has \textup{(SM)} if and only
if $X\times\R$ has \textup{(SM)}.
\end{theorem}

\begin{proof}
The reverse implication follows from the preceding observation.
For the forward implication, the case $X=\{0\}$ follows from
Theorem~\ref{thm:normed}. Suppose that $X\ne\{0\}$ has (SM).
Let $C\subset X\times\R$ be nonempty, closed, and convex, and fix
$z_0\notin C$. Strict separation gives
$0\ne L=(u^*,\alpha)\in X^*\times\R$ and $M\in\R$ such that
\begin{equation}\label{eq:product-separation}
 L(c)\le M<L(z_0)\qquad(c\in C).
\end{equation}
Choose $e=(v,b)$ with $b\ne0$ and $L(e)=-1$. If $u^*=0$, take
$e=(0,-1/\alpha)$; otherwise take $b=1$ and choose $v$ with
$u^*(v)=-1-\alpha$.

Set $D=C+\R_+e$, where $\R_+=[0,\infty)$. Then $D$ is convex and
$L(d)\le M$ for every $d\in D$, so $z_0\notin D$.
To prove closedness, take a net $z_i=c_i+s_i e\to z$, with $c_i\in C$
and $s_i\ge0$. Since
\[
 0\le s_i=L(c_i)-L(z_i)\le M-L(z_i),
\]
the net $(s_i)$ is eventually bounded and has a convergent subnet,
$s_{i_j}\to s\ge0$. Then
$c_{i_j}=z_{i_j}-s_{i_j}e\to z-se\in C$, and hence $z\in D$.
Thus $D$ is closed.

Use the continuous linear isomorphism
\[
 T(x,t)=(x+tv,bt),\qquad
 T^{-1}(y,s)=\bigl(y-(s/b)v,s/b\bigr).
\]
It satisfies $T(0,1)=e$ and $L(T(x,t))=u^*(x)-t$.
Therefore $A=T^{-1}(D)$ is nonempty, closed, and convex, and
\[
 A+\R_+(0,1)\subset A,\qquad
 (x,t)\in A\ \Longrightarrow\ t\ge u^*(x)-M.
\]
Each nonempty vertical section $A_x=\{t\mid(x,t)\in A\}$ is thus a
closed upper ray bounded below. Defining
\[
 \varphi(x)=\inf A_x\quad\text{with}\quad\inf\varnothing=+\infty,
\]
we obtain $A=\Epi\varphi$ and $\varphi\in\Gamma_\tau(X)$.
By Theorem~\ref{ch}(f), $A$ is portable.

Portability is preserved by continuous linear isomorphisms, because
supporting half-spaces and their contact points are transported by the
isomorphism. Hence $D=T(A)$ is portable. Since $z_0\notin D=D^\#$,
there are $d=c+se\in D$, with $c\in C$, $s\ge0$, and $n\in N_D(d)$
such that $\pair{z_0-d}{n}>0$.

The points $d+te$, $t\ge0$, belong to $D$, so $\pair{e}{n}\le0$.
Also $c=d-se\in D$, giving $-s\pair{e}{n}\le0$. Thus
$s\pair{e}{n}=0$. For every $y\in C\subset D$,
\[
 \pair{y-c}{n}=\pair{y-d}{n}+s\pair{e}{n}\le0,
\]
so $n\in N_C(c)$. Moreover,
\[
 \pair{z_0-c}{n}=\pair{z_0-d}{n}+s\pair{e}{n}>0.
\]
Therefore $z_0\notin C^\#$. Since this holds for every $z_0\notin C$,
we have $C=C^\#$. Theorem~\ref{ch} now gives (SM) on $X\times\R$.
\end{proof}

\begin{theorem}\label{cor:finite-products}
Let $X$ be a real Hausdorff locally convex space and $Y$ a
finite-dimensional real Hausdorff topological vector space. Then $X$
has \textup{(SM)} if and only if $X\times Y$ has \textup{(SM)}.
\end{theorem}

\begin{proof}
Write $n=\dim Y$. Since $Y$ is linearly homeomorphic to $\R^n$ and
(SM) is invariant under linear homeomorphisms, the forward implication
follows by applying Theorem~\ref{XxR} successively to the $n$ real
factors. The case $n=0$ is immediate. The reverse implication follows
from the factor argument at the beginning of this section.
\end{proof}

\begin{remark}
The argument above does not settle the product question for two arbitrary
infinite-dimensional (SM) spaces. By contrast, arbitrary countable
products do not preserve (SM): every factor $\R$ has (SM), but
$\omega=\R^\N$ fails it, as witnessed by \eqref{eq:omega-example}.
\end{remark}

\end{document}